\documentclass[10pt,a4paper]{article}

\usepackage{helvet}
\usepackage{amsmath,amsxtra,amssymb,xr}

\newtheorem{theorem}{Theorem}[section]
\newtheorem{lemma}[theorem]{Lemma}
\newtheorem{prop}[theorem]{Proposition}
\newtheorem{corollary}[theorem]{Corollary}

\newcommand{\GL}{\operatorname{GL}}
\newcommand{\Ind}{\operatorname{Ind}}
\newcommand{\Gal}{\operatorname{Gal}}
\newcommand{\Res}{\operatorname{Res}}
\newcommand{\Tr}{\operatorname{Tr}}
\newcommand{\Q}{\mathbb{Q}}
\newcommand{\Qbar}{\overline{\mathbb{Q}}}
\newcommand{\QFT}{\mathbb{Q}_{FT}}
\newcommand{\GFT}{G_{FT}}

\newcommand{\C}{\mathbb{C}}
\newcommand{\Z}{\mathbb{Z}}
\newcommand{\p}{\mathfrak{p}}
\newcommand{\PP}{\mathfrak{P}}
\newcommand{\CC}{\mathfrak{C}}
\newcommand{\Lp}{\mathbf{L}_p}
\newcommand{\Lptilde}{\widetilde{\mathbf{L}}_p}
\newcommand{\Lcal}{\mathcal{L}}

\newcommand{\Ad}{\mathrm{Ad}}
\newcommand{\OCp}{\mathcal{O}_{\mathbb{C}_p}}
\newcommand{\MCp}{\mathfrak{M}_{\mathbb{C}_p}}

\newcommand{\f}{\mathbf{f}}
\newcommand{\g}{\mathbf{g}}
\newcommand{\ord}{\mathrm{ord}}
\newcommand{\unit}{( \,p- \! \mathrm{adic \; unit})}
\newcommand{\period}{\Omega^+_E \, \Omega^-_E}

\newcommand{\dd}{\mathfrak{d}}
\newcommand{\qq}{\mathfrak{q}}
\newcommand{\rr}{\mathfrak{r}}

\newcommand{\OO}{\mathcal{O}}
\newcommand{\cc}{\mathfrak{c}}
\newcommand{\mm}{\mathfrak{m}}
\newcommand{\nn}{\mathfrak{n}}

\newcommand{\MM}{\mathfrak{M}}

\newcommand{\FF}{\mathbf{F}}
\newcommand{\G}{\mathbf{G}}
\newcommand{\Gcal}{\mathcal{G}}
\newcommand{\h}{\mathbf{h}}

\newcommand{\aaa}{\mathfrak{a}}
\newcommand{\bb}{\mathfrak{b}}

\begin{document}

\title{Congruences for Convolutions of Hilbert Modular Forms}
\date{}

\author{Thomas Ward}

\maketitle

\begin{abstract}
Let $\f$ be a primitive, cuspidal Hilbert modular form of parallel weight. We investigate the Rankin convolution $L$-values $L(\f,\g,s)$, where $\g$ is a theta-lift modular form corresponding to a finite-order character. We prove weak forms of Kato's `false Tate curve' congruences for these values, of the form predicted by conjectures in non-commmutative Iwasawa theory.
\end{abstract}
\section{Introduction}\label{Intro}
In recent years there has been much interest in the generalisation of Iwasawa theory to non-abelian field extensions. Let $p$ be an odd prime, $E$ an elliptic curve defined over $\Q$, and $F_\infty/\Q$ a $p$-adic Lie extension. In the paper \cite{mainconjecture}, Coates et al conjecture the existence of a non-abelian $p$-adic $L$-function in $K_1(\Z_p[[G]]_\mathcal{S^*})$ which interpolates the twisted $L$-functions $L(E,\rho,s)$ at $s=1$ (modified by certain simple factors). Here, $\rho$ ranges over the set of Artin representations of $G = \Gal(F_\infty/\Q)$, and $\Z_p[[G]]_\mathcal{S^*}$ is the localisation of $\Z_p[[G]]$ at a certain Ore set $\mathcal{S}^*$. More general conjectures of this nature were made by Fukaya and Kato in \cite{FukayaKato}.
\par
Attacking this appears very difficult in general. However, in the case of the `false Tate curve' extension of $\Q$, Kato proves in \cite{Kato} that the existence of the non-abelian $p$-adic $L$-function is equivalent to a set of strong congruences between certain abelian $p$-adic $L$-functions. Investigating these conjectural congruences is the motivation for our paper, and we will explain them further below.
\par
The false Tate curve extension is defined by
\[   
\QFT := \bigcup_{n \geq 1} \Q\left(\mu_{p^{n}},\sqrt[p^n]{\Delta}\right)  
\]
where $\mu_{p^n}$ denotes the group of $p^n$-th roots of unity, and $\Delta$ is a $p$-power free integer. The Galois group $\GFT := \Gal(\QFT/\Q)$ is a semi-direct product of two $p$-adic Lie groups of dimension one:
\begin{eqnarray*}
\GFT \; \cong \;  \left( \begin{array}{ccc} \Z_p^\times & 
\Z_p \\ 0 & 1 \end{array} \right) \; \lhd \; \GL_2(\Z_p) .
\end{eqnarray*}
This group has a unique self-dual representation of dimension $p^n-p^{n-1}$ (as is shown in \cite{VladThanasis} for example) which we denote by $\rho_{n,\Q}$. Putting $K_n= \Q(\mu_{p^n})$, this may be written 
\[
\rho_{n,\Q} \;=\; \Ind_{K_n}^{\Q} \chi_n
\]
for a one-dimensional character $\chi_n$ of $\Gal(K_n(\sqrt[p^n]{\Delta})/K_n)$.
In fact, all irreducible representations of $\GFT$ have the form $\rho_{n,\Q} \otimes \psi$ for some $n\geq 0$ and some finite-order character $\psi$ of $U^{(n)}$, where $U^{(n)} \;\cong\; \ker ( \Z_p^\times \twoheadrightarrow (\Z/p^n\Z)^\times ).$
\par
The structure of these Artin representations allows us to use the theory of Hilbert modular forms to make further progress. Let us write $F_n$ for the maximal real subfield of $K_n$. We have a Hilbert modular form over $F_n$ obtained as a theta-lift of $\chi_n$ (as defined in \cite[\S5]{Shimura}) which we denote by $\g_{\rho_n}$, identifying it with the two-dimensional induced representation $\rho_n := \Ind_{K_n}^{F_n} \chi_n$. Additionally, we know by the work of Wiles et al that $E/\Q$ is modular, and has an associated cusp form $f_E$. If we write $\f_E$ for the automorphic base-change of $f_E$ to the field $F_n$, then the non-abelian twist $L(E,\rho_{n,\Q},s)$ is essentially equal to the Rankin convolution $L(\f_E,\g_{\rho_n},s)$.
\par
This approach was used by Bouganis and V. Dokchitser in \cite{VladThanasis} to prove algebraicity properties for these $L$-values. It was then used in \cite{DelbourgoWard} by Delbourgo and the author to construct an abelian $p$-adic $L$-function $\Lp(E,\rho_n) \in \Z_p[[U^{(n)}]]$ interpolating the values $L(E,\rho_{n,\Q} \otimes \psi,1)$ for characters $\psi : U^{(n)} \rightarrow \Qbar^\times$.
\par
In this case, Kato proved in \cite{Kato} that the existence of the non-abelian $p$-adic $L$-function is equivalent to a family of congruences between the elements $\Lp(E,\rho_n)$. To be precise, there exists a map
\[   
\Theta_{G, \mathcal{S}^*} :  K_1 (\Z_p [[G]]_{\mathcal{S}^*}) \; \longrightarrow \; \prod_{n\geq 0} \mathrm{Quot}(\Z_p[[U^{(n)}]])^\times    
\]
whose image contains a sequence $(a_n)_{n\geq 0}$ if and only if 
\[     
\prod_{1\leq j \leq n} N_{j,n}\left(  \frac{a_j}{N_{0,j}(a_0)}. \frac{\varphi \circ N_{0,j-1}(a_0)}{\varphi(a_{j-1})} \right)^{p^j}  \, \equiv \, 1 \mod p^{2n}  \qquad \text{ for all }  n \in \mathbb{N} . 
\]
Here we must explain that $N_{i,j} : \Z_p[[U^{(i)}]]^\times\rightarrow \Z_p[[U^{(j)}]]^\times $ denotes the norm map, and $\varphi : \Z_p[[\Z_p^\times]] \rightarrow \Z_p[[\Z_p^\times]]$ is the ring homomorphism induced by the $p$-power map on $\Z_p^\times$. 
\par
In the case $p=3$, the first level of these congruences for $a_n = \Lp(E,\rho_n)$ has been verified by Bouganis in \cite{BouganisJLMS}. In particular, he uses deep results of Wiles to prove a relationship between the motivic period $\period$ appearing in the conjectures, and an automorphic period associated to the modular form $f_E$.
\par
Delbourgo and the author also proved a congruence of this form (using Hilbert modular forms as described above) but modulo a smaller power of $p$. We did this for semistable elliptic curves in \cite{DelbourgoWard}, and extended our results to the case of CM curves in \cite{DelbourgoWard2}.
\par
In this paper, we will show that one may take a primitive Hilbert modular form $\f$ of arbitrary parallel weight, and achieve similar results for the convolution $L(\f,\g_{\rho_n},s)$ at all critical values. 
Let $\f$ be a primitive, cuspidal Hilbert modular form over $F_n$, with parallel weight $k \geq 2$, conductor $\cc(\f)$ and Hecke character $\eta$.
We write $\p$ for the unique prime ideal of $\OO_{F_n}$ above $p$. 
To state our main results, we must impose the following hypotheses.
\\\\
\textit{Hypothesis (Ord):} $\p$ does not divide $\cc(\f)$ or $\Delta \OO_F$, and $\cc(\f)+ \Delta \OO_F = \OO_F$. Further, the Fourier coefficient $C(\cc(\f),\f)$ is non-zero and $C(\p,\f)$ is a $p$-adic unit.
\\\\
{\it Hypothesis (Cong):} there exists no congruence modulo $\MCp$ between $\f$ and another Hilbert modular form which lies outside the $\f$-isotypic component of $\mathcal{M}_k(\cc(\f),\eta)$.
Here, $\MCp$ denotes the maximal ideal of $\OCp$.
\begin{theorem}\label{theorem1}
Suppose that $p>k-2$, that $\f$ has rational Fourier coefficients and that Hypotheses (Ord) and (Cong) are satisfied. For each critical value $1 \leq r \leq k-1$ there exists a unique element $\Lp(\f,\rho_n,r) \in \Z_p[[U^{(n)}]]$ with the property
\[
\psi \Big(  \Lp(\f,\rho_n,r)  \Big) 
\;=\;  
\frac{\epsilon_{F_n}(\rho_n\otimes\psi,1-r)_\p}{ \alpha(\p)^{f(\rho_n\otimes\psi, \p)}} 
\;\times \;
\frac{ P_\p(\rho_n\otimes\psi^{-1}, \, \alpha(\p)^{-1} p^{r-1}) }{ P_\p(\rho_n\otimes\psi, \, \alpha(\p) \, p^{-r})} 
\]
\[
\;\times \;
\frac{ \Psi_S(\f, \g_{\rho_n\otimes\psi}^\iota, r) }{ D_F^{k-2} \; \left< \f,\f \right>_{\cc(\f)} } 
\]
for each character $\psi : U^{(n)}\rightarrow \C^\times$. 
\end{theorem}
Here we write $\Psi_S(\f, \g_{\rho_n\otimes\psi}^\iota, r)$ for the completed Rankin convolution of $\f$ and $\g_{\rho_n\otimes\psi}^\iota$ with the Euler factors at primes dividing $p.\Delta\OO_F$ removed, and $\left< \f,\f \right>_{\cc(\f)}$ for the Petersson self-product of $\f$. The other terms in the formula are defined in \S\ref{section2} and \S\ref{section6}.
\begin{theorem}\label{theorem2}
Suppose that $\f$ satisfies the same hypotheses as in Theorem \ref{theorem1}. 
Put $a_j= \Lp(\f,\rho_j,r)$ for some $1 \leq r \leq k-1$. Then we have the congruence
\[     
\prod_{1\leq j \leq n} N_{j,n}\left(  \frac{a_j}{N_{0,j}(a_0)}. \frac{\varphi \circ N_{0,j-1}(a_0)}{\varphi(a_{j-1})} \right)^{p^j}  \, \equiv \, 1 \mod p^{n+1}  
\]
for each $n\geq 0$.
\end{theorem}
In Theorems \ref{theorem1} and \ref{theorem2}, we make the assumption $p>k-2$. This is forced on us for technical reasons, and we expect it is not neccessary for the results to hold. In \S\ref{section7} we will give two numerical examples for $p=3$, in which the first level congruences from Theorem \ref{theorem2} hold even without this condition.
\par
However, we will also present an example in which Hypothesis (Cong) fails  and the congruences fail. This is not surprising, given a well-known connection between our complex period $\left< \f,\f \right>_{\cc(\f)}$ and the congruence primes of $\f$ (we discuss this at the end of \S\ref{section6}).
\\\\
\noindent
\textit{Acknowledgements:}
The author thanks Daniel Delbourgo and Antonio Lei for many useful suggestions.
\section{Convolutions of Hilbert Modular Forms}\label{section2} 
In \S\ref{section2} and \S\ref{section3} we review some results on convolutions of Hilbert modular forms, and the Rankin-Selberg method. In \S\ref{section4} we will go on to construct $p$-adic measures which interpolate the critical values of these $L$-functions. Our principal reference for this theory is Panchishkin's book \cite{Pan}, and our methods are heavily based on his.
\par
Let $p$ be an odd prime, and let $F$ be a totally real number field. We fix a prime ideal $\p$ of $F$ which lies above $p$.
Consider two Hilbert modular forms $\f$ and $\g$ defined over $F$, which are primitive and cuspidal. Suppose that $\f$ has parallel weight $k \geq 2$, conductor $\cc(\f)$ and Hecke character $\eta$; also suppose that $\g$ has parallel weight one, conductor $\cc(\g)$ and character $\omega$. Weights of Hilbert modular forms will always be assumed parallel in what follows.
\par
We will adopt the following hypotheses throughout:
\[     \p \nmid \cc(\f), \quad \cc(\f)+\cc(\g) = \OO_F , \quad C\big(\cc(\f),\f\big) \neq 0   \]
where we have written $C(\mm, \f)$ for the Fourier coefficient of $\f$ at the ideal $\mm$. Further, we adopt the hypothesis that $\f$ is $\p$-ordinary: by this we mean that $C(\p,\f)$ is a $p$-adic unit. 
\par
The $L$-series associated to $\f$ is defined by
\[      
L(\f,s)   \; := \;    \sum_\mm  \, C(\mm,\f)    \, N(\mm)^{-s} .     
\]
As we assume $\f$ is primitive, we also have the Euler product expression:
\[      
L(\f,s)   \; = \;   \prod_\qq \left( 1 \, - \,  C(\qq,\f) \, N(\qq)^{-s}  \, + \,  \eta(\qq) \, N(\qq)^{k-1-2s} \right)^{-1} . 
\]
We write $\Psi(\f,\g,s)$ for the completed Rankin convolution of $\f$ and $\g$, which is given by
\[  \Psi(\f,\g_\rho,s) \, := \, \left( \frac{\Gamma(s)}{(2\pi)^{ s }}  \right)^{2[F:\Q]}   \,   L_{\cc(\f)\cc(\g)}(2s-k-1, \eta \omega) \, \sum_{\aaa} C(\aaa, \f) \, C(\aaa, \g)\, N(\aaa)^{-s}  . 
\]
These $L$-series only converge for $\mathrm{Re}(s)$ sufficiently large, but both may be continued to holomorphic functions on the whole complex plane, and satisfy functional equations of the usual form (see \cite{Shimura} for example).
\par
For an integral ideal $\aaa$ of $\OO_F$ we have two linear operators $|\aaa$ and $|U(\aaa)$ on the space $\mathcal{M}_k(\cc(\f), \eta)$, which may be defined by their effect on the Fourier coefficients of any Hilbert modular form $\h$:
\[ 
      C(\mm, \h \big| \aaa)  \; = \;  C(\mm \aaa^{-1}, \h )  \qquad \mathrm{and} \qquad  C(\mm, \h| U(\aaa)) = C(\mm \aaa, \h )  , 
\]
where we put $C(\mm, \h)=0$ whenever the ideal $\mm$ is not integral. Here we also have an involution $J_\aaa$ on this space, defined in \cite[Chapter 4]{Pan}.
\begin{theorem}\label{rep1}
Put $d=[F:\Q]$, and let $\dd$ be the different of $F/\Q$. Then, for any $\FF \in \mathcal{S}_k(\cc(\FF),\eta)$ and $\G \in \mathcal{S}_1(\cc(\G),\omega)$, we have the following integral representation for the Rankin convolution: 
\[
\Psi(\FF,\G,s)  \;\; = \;\;  (-1)^{d(s-k+1)}   \;   2^{dk}   \;   i^{d(1-k)}   \;  N\big( \cc(\FF) \cc(\G) \dd^2 \big)^{(k-1)/2-s}   \;  N\big( \cc(\G) \big)^{1-k/2} 
\]
\[
\times \quad\left< \FF^\iota , \Big( \G| J_{\cc(\FF) \cc(\G)} \,.\, E_{k-1} (s-k+1, \eta \omega^{-1}) \Big) \Big| U\big(\cc(\G)\big) \circ J_{\cc(\FF)}  \right>_{\cc(\FF)}   ,
\]
where $E_{k-1}(s,\psi)$ is the Eisenstein series specified in \cite[chapter 4]{Pan}.
\end{theorem}
{\it Proof.} 
We take the integral representation from \cite[4.32]{Shimura}, then apply the trace operator $\Tr^{\cc(\FF)\cc(\G)}_{\cc(\FF)}$. Using the identity
\[
V \big|\Tr^{\aaa \bb}_{\aaa} \; = \; N(\bb)^{1-k/2}  \; \, V \big| J_{\aaa\bb} \circ U(\bb) \circ J_{\aaa} ;
\]
we obtain the desired result, following the same calculation as in \cite[page 136, \S4.4]{Pan}.
\par
We return to our primitive Hilbert modular forms $\f$ and $\g$. For a prime ideal $\qq$ of $F$ we will always write
 \[ 1 - C(\qq,\f)X + \eta(\qq) N(\qq)^{k-1} X^2  \; = \;  \big( 1-\alpha(\qq)X \big) \big( 1-\alpha'(\qq)X \big) , \] 
for the factorisation of the local polynomial of $\f$ at $v$. For the prime $\p$ we choose $\alpha(\p)$ to be the root which is a $p$-adic unit, with $\alpha'(\p)$ the non-unit root (which we can do as we assumed $\f$ to be ordinary at $\p$). Similarly, for $\g$ we will write
\begin{eqnarray*} 
1 - C(\qq,\g)X + \omega (\qq) X^2  &=&  \big(1-\beta(\qq)X\big) \big(1-\beta'(\qq)X\big), 
\\ 
1 - \overline{C(\qq,\g)}X + \omega^{-1} (\qq) X^2  &=& \big(1-\hat{\beta}(\qq)X\big)  \big(1-\hat{\beta}'(\qq)X\big) . 
\end{eqnarray*}
The convolution of $\f$ and $\g$ may be written as the following Euler Product:
\[
L_{\cc(\f)\cc(\g)}(2s-k-1, \eta \omega) \, \sum_{\aaa} C(\aaa, \f) \, C(\aaa, \g)\, N(\aaa)^{-s}  
\;\,=\;\,
\prod_{\qq}  (1- \alpha(\qq) \beta(\qq) N(\qq)^{-s})^{-1} 
\]
\[
\times \;\; (1- \alpha(\qq) \beta'(\qq) N(\qq)^{-s})^{-1} \; (1- \alpha'(\qq) \beta(\qq) N(\qq)^{-s})^{-1} \; (1- \alpha'(\qq) \beta'(\qq) N(\qq)^{-s})^{-1} .
\]
Let us now fix a squarefree ideal $\mm_0$ which is divisible by $\p$ and all the primes dividing $\cc(\g)$. We define the \mbox{$\mm_0$-stabilisation} of $\f$ to be
\[
     \f_0 \;:=\; \sum_{\aaa| \mm_0} \mu(\aaa) \, \alpha'(\aaa) \,.\, \f \big| \aaa      
\]
where $\mu$ is the M\"obius function on ideals. This definition is equivalent to the identity
\[
      L(\f_0,s)   \; = \;    L(\f,s)  \; \times \;  \prod_{\qq|\mm_0}   (1-\alpha'(\qq) N(\qq)^{-s})    .    
\]
We also define $\g_{\mm_0} \in \mathcal{S}(\cc(\g)\mm_0^2,\omega)$ by
\[
\g_{\mm_0}        \; = \;    \sum_{\aaa|\mm_0} \mu(\aaa) \cdot \, \g \big| U(\aaa) \circ \aaa .
\]
Equivalently, $\g_{\mm_0}$ is the non-primitive Hilbert modular form whose Fourier coefficients are given by
\begin{equation*}
 C(\nn,  \g_{\mm_0}) = 
\begin{cases}
  C(\nn,  \g)           & \text{ if $\nn$ and $\mm_0$ are coprime;} \\
                     0          & \text{ otherwise.}
\end{cases}
\end{equation*}
For the rest of this section, we will write $\mm'$ for an auxiliary ideal supported on the primes dividing $\mm_0$, such that $\cc(\g)\mm_0^2|\mm'$. We substitute
\[
      \FF \;=\; \f_0 \;\in\; \mathcal{S}_k \big(\cc(\f)\mm_0,\eta \big)  , \; \quad \mbox{and}  \; \quad   \G \;=\; \g_{\mm_0} |J_{\cc(\f) \mm'} \;\in\; \mathcal{S}_1 \big( \cc(\f) \mm', \omega^{-1} \big)       
\]
into Theorem \ref{rep1} to obtain the formula
\[
       \Psi\big( \f_0, \g_{\mm_0} |J_{\cc(\f) \mm'}, s\big)  \;=\;  (-1)^{dk}   \; 2^{dk} \, i^{d(1-k)} \; N(\mm'\mm_0^{-1})^{1-k/2}    \;   N\big( \cc(\f)  \mm' \dd^2 \big)^{(k-1)/2-s} 
\]
\[
    \times  \quad \left< \f_0^\iota , \; \g_{\mm_0} \,.\, E_{k-1} (s-k+1, \eta \omega^{-1}) \Big| U(\mm'\mm_0^{-1}) \circ J_{\cc(\f)\mm_0}  \right>_{\cc(\f)\mm_0} .    
\]
We define a linear functional
\begin{eqnarray*}   
\Lcal_F   \;:\; \mathcal{M}_k (\cc\mm_0, \eta) &\longrightarrow& \C  \\
                  \Phi  &\longmapsto&  \frac{ \left< \f_0^\iota , \Phi|J_{\cc(\f) \mm_0} \right>_{\cc(\f) \mm_0} }{ \left< \f,\f \right>_{\cc(\f)} }  .
\end{eqnarray*}           
We also write $\mathcal{H}ol$ for the holomorphic projection operator constructed in \cite[page 138, \S4.6]{Pan}. 
This operator maps the space $\tilde{\mathcal{M}}_k(\cc\mm_0,\eta)$ of $C^\infty$-Hilbert modular forms to $\mathcal{M}_k(\cc\mm_0,\eta)$, and is related to the Petersson inner product by the formula
\[
\left< \h , \Phi \right>_{\cc\mm_0}  \;=\;  \left< \h , \mathcal{H}ol (\Phi) \right>_{\cc\mm_0}
\]
for any $\h$ in $\mathcal{S}_k(\cc\mm_0,\eta)$.
Holomorphic projection was not required in \cite{DelbourgoWard}, where the only weight considered was $k=2$; however the Eisenstein series $E_{k-1}(s-k+1)$ will be non-holomorphic at some critical values if $k>2$. For convenience we will put
\[   \Phi(\g, s)  \;:=\;  \mathcal{H}ol \left( \g_{\mm_0} E_{k-1}(s,\eta\omega^{-1}) \right)  . \]
Then, applying holomorphic projection to the above formula for $\Psi( \f_0, \g_{\mm_0} |J_{\cc(\f) \mm'}, s)$ and writing it in terms of the linear functional $\Lcal_F$, we obtain
\[    \frac{ \Psi(\f_0, \g_{\mm_0} |J_{\cc\mm'}, s) }{ \left< \f,\f \right>_{\cc(\f)} } \;=\; (-1)^{dk}   \; 2^{dk} i^{d(1-k)} \; N(\mm'\mm_0^{-1})^{1-k/2}  
   N(\cc(\f)\mm' \dd^2)^{(k-1)/2-s} \] 
\[    \times \;\;\;\; \Lcal_F \left( \Phi( \g, s-k+1) \big| U(\mm'\mm_0^{-1}) \right) .   \]
Next we will rewrite the convolution $\Psi\big( \f_0, \g_{\mm_0} |J_{\cc(\f) \mm'}, s\big)$ in terms of $\Psi(\f,\g^\iota,s)$. We define the contragredient Euler factor by 
\begin{eqnarray*}
  \mathrm{Eul}_{\mm_0}(\g^\iota, s) &:=& \prod_{v | \mm_0}  (1-\alpha'(v)\hat{\beta}(v)N(v)^{-s}) 
(1-\alpha'(v)\hat{\beta}'(v)N(v)^{-s})  \\
 && 
 \quad \times \; (1-\alpha^{-1}(v)\beta(v)N(v)^{s-1}) (1-\alpha^{-1}(v)\beta'(v)N(v)^{s-1}) .
\end{eqnarray*}
\begin{lemma}\label{lemma2}
We have the formula
\begin{eqnarray*}    \Psi(\f_0, \g_{\mm_0}|J_{\cc\mm'}, s)    &=&    N \left( \frac{\cc(\f)\mm'}{\cc(\g)} \right)^{1/2-s}   \frac{ \alpha(\mm') }{ \alpha\big(\cc(\g)\big) } \;
\\
  & & \times \;\;\;\;\;\;  \Lambda(\g) \;\;  C\big(\cc(\f),\f\big) \;\; \mathrm{Eul}_{\mm_0}(\g^\iota,s) \;\; \Psi(\f, \g^\iota, s)  . 
\end{eqnarray*}
\end{lemma}
{\it Proof.} 
Let us put $\FF=\f_0$ and $\G=\g_{\mm_0}|J_{\cc\mm'}$. 
Quoting \cite[page 125, 1.22]{Pan} we have the identity
\[
      \h \big| J_{\mm\cc(\h)}  \;=\;  N(\mm)^{\mathrm{wt}(\h)/2}  \left( \h \big|J_{\cc(\h)} \right) \big|\mm   
\]
which holds for all Hilbert modular forms $\h$ and ideals $\mm$. 
\par
By assumption, $\mm'= \cc(\g) \mm_0^2 \rr$ for some integral ideal $\rr$ which is supported on the primes dividing $\mm_0$. The level of $ \g_{\mm_0}$ is $\cc(\g)\mm_0^2$, so we apply the above identity to obtain
\[  \G \; = \; \g_{\mm_0} \big| J_{\cc(\f) \cc(\g) \mm_0^2 \rr}   
\;=\;   
 N \big(  \cc(\f) \rr  \big)^{1/2}    \left( \g_{\mm_0} \big| J_{\cc(\g)\mm_0^2} \right)   \Big|  \cc(\f) \rr     \]
Let us write $\tilde{\g} = \g_{\mm_0} \big| J_{\cc(\g)\mm_0^2}$ so that $\G = N (\cc(\f) \rr )^{1/2} \;\tilde{\g} \,|\cc(\f)\rr$. Then we have
\begin{eqnarray*}    
\Psi(\FF, \G, s)   &=&   N \big( \cc(\f)\rr \big)^{1/2} \Psi \big(\f_0, \tilde{\g} |\cc(\f)\rr, s \big)    \\
     &=&  
     N \big( \cc(\f)\rr \big)^{1/2-s}   \Psi \big( \f_0 |U( \cc(\f)\rr ), \tilde{\g}, s \big)      \\
     &=&  
     N \big( \cc(\f)\rr \big)^{1/2-s}  \alpha( \rr ) C\big(\cc(\f),\f\big)  \Psi \big(\f_0, \tilde{\g}, s \big)     .
\end{eqnarray*}
The second equality here follows from the identity
\[
\Psi \big(\f,\g |\aaa,s \big) \;=\; N(\aaa)^{-s} \, \Psi \big( \f |U(\aaa),\g,s \big)
\]
which is clear when we recall that $C(\mm, \g | \aaa)  \; = \;  C(\mm \aaa^{-1}, \g)$ and $C(\mm, \f| U(\aaa)) = C(\mm \aaa, \f )$. 
The third equality is deduced from the fact that $\f_0|U(\qq) = \alpha(\qq) \,\f_0$ when $\qq$ divides $\mm_0$ or $\cc(\f)$,
as well as the observation $\alpha(\cc(\f))=C(\cc(\f),\f)$.
We also have
\begin{eqnarray*}    
  \Psi(\f_0, \tilde{\g}, s)   &=&   \Psi(\f_0, \g_{\mm_0}|J_{\cc(\g)\mm_0^2}, s)   \\
     &=&  
   N(\mm_0)^{1-2s} \;  \alpha(\mm_0)^2 \;  \Lambda(\g) \; \mathrm{Eul}_{\mm_0}(\g^\iota,s) \; \Psi(\f, \g^\iota, s) ,     
\end{eqnarray*}
as quoted from \cite[page 130, proposition 3.5]{Pan}. Combining these two equations, we get the desired result.
\par
Combining Lemma \ref{lemma2} with our integral representation, we obtain 
\begin{eqnarray*}    
  \frac{ \Psi\big( \f_0, \g_{\mm_0} |J_{\cc(\f) \mm'}, s\big) }{ \left< \f,\f \right>}  &=&  (-1)^{dk}   \; 2^{dk} i^{d(1-k)} \; N(\mm'\mm_0^{-1})^{1-k/2}    \;   N(\cc(\f)\mm' \dd^2)^{(k-1)/2-s}  \\ 
  &&    \times  \qquad   \Lcal_F \left( \Phi( \g, s-k+1) \big| U(\mm'\mm_0^{-1}) \right)   
\end{eqnarray*}
\[
= \;\;\;\; N\left( \frac{\cc(\f)\mm'}{\cc(\g)} \right)^{1/2-s} \; \Lambda(\g) \; \frac{\alpha(\mm')}{\alpha\big(\cc(\g)\big)}  \;  C\big(\cc(\f),\f\big)  \; \mathrm{Eul}_{\mm_0}(\g^\iota,s) 
        \; \frac{ \Psi(\f, \g^\iota, s) }{ \left< \f,\f \right>} 
\]
This implies
\[ 
\alpha\big(\cc(\g)\big)^{-1} \; \Lambda(\g) \; \mathrm{Eul}_{\mm_0}(\g^\iota,s) \; \frac{ \Psi(\f, \g^\iota, s) }{ \left< \f,\f \right>_{\cc(\f)} }  
 \;=\;  
  \frac{(-1)^{dk} 2^{dk} i^{d(1-k)}}{\alpha(\mm') \, C\big(\cc(\f),\f\big)}   N(\mm'\mm_0^{-1})^{1-k/2}   
\]
\[
   \times \;\;  N\big( \cc(\f)\mm'\dd^2 \big)^{(k-1)/2-s} \;  N\left( \frac{\cc(\f)\mm'}{\cc(\g)} \right)^{s-1/2}    \;\;  \Lcal_F \Big( \Phi( \g, s-k+1) \big| U(\mm'\mm_0^{-1}) \Big)  
\]
and we get the following corollary.
\begin{corollary}\label{corol3}
We have the formula
\[   
 N\big( \cc(\g)\dd^2 \big)^{s-1/2}  \;\; \Lambda(\g) \;\;  \alpha\big(\cc(\g)\big)^{-1}  \;\; \mathrm{Eul}_{\mm_0}(\g^\iota,s) 
 \;\;
 \frac{ \Psi(\f, \g^\iota, s) }{ \left< \f,\f \right>_{\cc(\f)} } 
 \]
 \[
 =  \;\;\;
  \frac{ (-1)^{dk} \, 2^{dk} \, i^{d(1-k)} }{\alpha(\mm') \, C\big(\cc(\f),\f\big)} \;\;  N\big( \cc(\f)\mm_0 \dd^2 \big)^{k/2-1}  
     \;\;    \Lcal_F \Big( \Phi( \g, s-k+1) \big| U(\mm'\mm_0^{-1}) \Big)  .
\]
\end{corollary}
This is a more general version of \cite[corollary 11]{DelbourgoWard}.
We observe that the left-hand side of this equation does not depend on the choice of $\mm'$, so neither does the right-hand side; this is the `distribution property' which we will require in \S\ref{section4}. 
\section{Fourier Expansions}\label{section3}
We will need an explicit formula for the Fourier coefficients of the Hilbert modular form 
\[   \Phi( \g, s)  \;:=\;  \mathcal{H}ol \left( \g_{\mm_0} E_{k-1}(s,\eta\omega^{-1}) \right) .  \]
We now very briefly recall some facts about Fourier expansions of Hilbert modular forms (full details can be found in \cite[chapter 4]{Pan} or \cite{Shimura}). 
Let $h=|\tilde{Cl}_F|$ be the narrow ideal class number of the field $F$,
and choose finite ideles $t_1,...,t_h$ such that the set of ideals $\{ \tilde{t}_\lambda : 1 \leq \lambda \leq h \}$ form a complete set of representatives for $\tilde{Cl}_F$, and are all coprime to $\mm_0$. Here we write $\tilde{t}_\lambda$ for the ideal of $\OO_F$ generated by $t_\lambda$.
A Hilbert modular form $\f$ may be naturally identified with an $h$-tuple $(f_1,...,f_h)$ of modular forms on $\mathbb{H}^n$, where $\mathbb{H}$ denotes the complex upper half-plane.
For $1 \leq \lambda \leq h$ we will refer to $f_\lambda$ as the `$\lambda$-component' of $\f$. This function has a Fourier expansion of the form
\[
 f_\lambda (z) \;=  \sum_{\xi}  a_\lambda(\xi) \, e_F(\xi z)
\]
where the sum ranges over all $\xi \in \tilde{t}_\lambda$ which are totally positive (denoted `$0 \ll \xi$') and $\xi=0$.
Further, if the ideal $\mm = \xi \tilde{t}_\lambda^{-1}$ is integral, then $C(\mm,\f) = a_\lambda(\xi) N(\tilde{t}_\lambda)^{-k/2}$.
\par
Let us write the $\lambda$-component of $\Phi(\g,s)$ as
\[    \Phi( \g, s)_\lambda (z) \;=  \sum_{ 0 \,\ll\, \xi \in \,\tilde{t}_\lambda}  \phi_\lambda(\xi,\g,s) e_F(\xi z)   \]
and the $\lambda$-component of $\g$ as
\[    \g_\lambda (z) \;=  \sum_{ 0 \,\ll\, \xi \in \,\tilde{t}_\lambda}  b_\lambda(\xi) e_F(\xi z) .   \]
We quote \cite[page 143, 5.8]{Pan}, specialising to our case of $\mathrm{wt}(\g)=1$ as before: for an integer $1 \leq r \leq k-1$, we have
\begin{eqnarray*}    
 \phi_\lambda(\xi,\g,r-k+1)  &=&    N(\tilde{t}_\lambda)^{ \frac{k+1}{2}-1-r }  \sum_{\xi = \xi_1 + \xi_2} b_\lambda(\xi_1)  \;  \sum_{\tilde{\xi_2} \,=\, \tilde{b} \, \tilde{c}}
 \mathrm{sign}( N(\tilde{b}) )^{k-2}    \;      N(\tilde{b})^{2r-k} 
  \\
     &\times&   (\eta \omega^{-1})(\tilde{c})  \;  \prod_\nu   P_{r-k+1}(\xi_{2,\nu}, \xi_\nu)
\end{eqnarray*}
where the second sum ranges over all $b \in \tilde{t}_\lambda^{-1}$ and $c \in \OO_F$ such that $\tilde{\xi_2} = \tilde{b}\,\tilde{c}$.
Here, $\xi_{2,\nu}$ denotes the image of $\xi_2$ under the real embedding $\nu: F\hookrightarrow \mathbb{R}$, and $P_s(\xi_{2,\nu}, \xi_\nu)$ denotes the polynomial 
\[           P_{s}(\xi_{2,\nu}, \xi_\nu)   \;=\;  \sum_{i=0}^{-s}  (-1)^i \binom{-s}{i} \frac{ \Gamma(k-1+s) }{ \Gamma(k-1+s-i) }  
            \frac{ \Gamma(k-1-i) }{ \Gamma(k-1) }         \xi_{2,\nu}^{-s-i}  \; \xi_\nu^{\;i}  \]
where $s \leq 0$. For $s \in \Z$, this polynomial has rational coefficients. Now, recall that 
\begin{equation*}
 C(\mm,  \Phi( \g, s-k+1)) = 
\begin{cases}
   N(\tilde{t}_\lambda)^{-k/2} \; \phi_\lambda(\xi,\g,s-k+1)  & \text{ if $\mm=\xi \tilde{t}_\lambda^{-1}$ is integral;} \\
                                         0          & \text{ if  $\mm$ is not integral.}
\end{cases}
\end{equation*}
Therefore,
\begin{eqnarray*}    
 C(\mm,  \Phi(\g, s-k+1))   &=&    N(\tilde{t}_\lambda)^{ -\frac{1}{2}-s }  \sum_{\xi = \xi_1 + \xi_2} b_\lambda(\xi_1)  \;  \sum_{\tilde{\xi_2} = \tilde{b}\tilde{c}}
 \mathrm{sign}( N(\tilde{b}) )^{k}    \;      N(\tilde{b})^{2s-k} 
  \\
     &\times&   (\eta \omega^{-1})(\tilde{c})  \;  \prod_\nu   P_{s-k+1}(\xi_{2,\nu}, \xi_\nu)
\end{eqnarray*}
where $\mm=\xi \tilde{t}_\lambda^{-1}$ for $\xi \gg 0$ as above. If we substitute $\xi_2 = \xi-\xi_1$, we can write
\[  \prod_\nu P_{s}(\xi_{2,\nu}, \xi_\nu) \;=\; (-1)^{ds}\; N(\xi_1)^{-s} \;+ \; N(\xi) \,\times\, \mbox{(other terms)}  \]
and therefore
\[   
 C\big(\mm,  \Phi(\g, s-k+1)\big)   \,\;=\;\,    N(\tilde{t}_\lambda)^{ -\frac{1}{2}-s } \Big(   N(\xi) \, u_\lambda(\xi) 
\]
\[
 +\;   (-1)^{d(s-k+1)} \sum_{\xi = \xi_1 + \xi_2} b_\lambda(\xi_1) \; N(\xi_1)^{s-k+1} \sum_{\tilde{\xi_2} = \tilde{b}\tilde{c}}
 \mathrm{sign}( N(\tilde{b}) )^{k}    \;      N(\tilde{b})^{2s-k} \; (\eta \omega^{-1})(\tilde{c})       \Big)
\]
where $u_\lambda(\xi)$ is a linear combination of the coefficients $b_\lambda(\xi)$ and the values of the Hecke character $\eta\omega^{-1}$ (just as in \cite[page 143, 5.9]{Pan}).
\section{Bounded $p$-adic Measures}\label{section4}
Let us now fix an imaginary quadratic extension $K/F$, where $F$ is our totally real field as before. We will construct a bounded measure on the group 
\[
\Gcal(\MM_0) \; := \; \Gal \big( K(\MM_0 p^\infty)/K \big) 
\] 
where $\MM_0 = \mm_0 \OO_K$ and $K(\MM_0 p^\infty)$ denotes the maximal ray class field modulo $\MM_0 p^\infty$ over $K$.
\par
For every finite-order character $\chi : \Gcal(\MM_0) \rightarrow \C^\times$, we have a theta-lift Hilbert modular form $\g \in \mathcal{M}_1(\cc(\g),\omega)$ defined over $F$ (see \cite[\S5]{Shimura}). We denote this by $\g_\rho$, identifying it with the induced Artin representation 
\[
\rho = \Ind_{\Gal(M/K)}^{\Gal(M/F)} \chi
\]
where $M/K$ is a finite extension through which the character $\chi$ factors. This identification is natural since $L(\g_\rho,s)$ is the same as the Artin $L$-function $L(\rho/F,s)$. 
\par
In what follows we will write $\epsilon_F(\rho,s)$ for the global epsilon factor that appears in the functional equation of $L(\rho/F,s)$; to be precise we have
\[          
\hat{L}(\rho,s)  \; = \;  \epsilon_F(\rho,s)  \; \hat{L}(\rho^\vee,1-s)   
\]
where $\hat{L}(\rho,s)= L_\infty(\rho,s) \, L(\rho,s)$ is the $L$-function with the Euler factors at infinite places included, and $\rho^\vee$ is the contragredient representation. 
We will also write $\cc(\rho)$ for the Artin conductor of $\rho$, and
$\mathrm{Eul}_{\mm_0}(\rho^\vee,s)$ for the Euler factor $\mathrm{Eul}_{\mm_0}(\g_{\rho}^\iota, s)$ we defined in \S\ref{section2}.
\begin{prop}\label{measureprop}
Given a primitive cusp form $\f \in \mathcal{S}_k \big( \cc(\f),\eta \big)$ there exists an algebraic-valued, bounded measure $\mu_{\mm_0}(\f,r)$ (for each $1 \leq r \leq k-1$) on $\Gcal(\MM_0)$ taking the value
\[      \int_{\Gcal(\MM_0)} \chi \; d\mu_{\mm_0}(\f,r)  \;\;=\;\;    \epsilon_F(\rho,1-r) \; \mathrm{Eul}_{\mm_0}(\rho^\vee,r) \;   
        \alpha\big(\cc(\rho) \big)^{-1} 
\;\;  \frac{ \Psi(\f, \g_\rho^\iota, r) }{ \left< \f,\f \right>_{\cc(\f)} }    \]
at every finite-order character $\chi : \Gcal(\MM_0) \rightarrow \C^\times$, where $\rho = \Ind_K^F \chi$. 
\end{prop}
{\it Proof.} 
Firstly, we know that these values are algebraic, by results from the key paper \cite{Shimura} of Shimura.
It is a simple consequence of \cite[2.48]{Shimura} that
\[
\Lambda(\g_\rho) \; N(\cc(\g_\rho) \dd^2)^{r-1/2} \; i^{[F:\Q]}  \; = \;  \epsilon_F(\rho,1-r).
\]
Further, by definition of the cusp form $\g_\rho$, the conductor $\cc(\g_\rho)$ is equal to $\cc(\rho)$. Therefore, by Corollary \ref{corol3} we may write
\[   \epsilon_F(\rho,1-r) \;   \frac{ \mathrm{Eul}_{\mm_0}(\rho^\vee,r) }{ \alpha\big(\cc(\rho)\big) } 
\;  \frac{ \Psi(\f, \g_\rho^\iota, r) }{ \left< \f,\f \right>_{\cc(\f)} } \;=\; \gamma(\mm')  \;\, \Lcal_F \Big( \Phi( \g_\rho, r-k+1) \big| U(\mm'\mm_0^{-1}) \Big)     \]
for any $\mm'$ such that $\cc(\rho)\mm_0^2|\mm'$. Here we have written $\gamma(\mm')$ for the constant
\[      \gamma(\mm') \;=\; \frac{ (-1)^{d} \; 2^{dk} \; i^{dk} }{\alpha(\mm') \; C\big(\cc(\f),\f\big)}   N \big( \cc(\f)\mm_0\dd^2 \big)^{k/2-1}    . \]
Quoting \cite[page 144, 5.11]{Pan} we know that that the linear functional $\Lcal_F$ may be written as a linear combination of Fourier coefficients:
\[  \Lcal_F \left(  \Theta \right)  \;=\;   \sum_\aaa \kappa_\aaa \, C(\aaa, \Theta) ,  \]
where the $\kappa_\aaa$ are algebraic numbers, fixed independently of $\Theta$, and all but finitely many $\kappa_\aaa$ are zero. As Panchishkin states in \cite{Pan}, this follows from a version of Atkin-Lehner theory for Hilbert modular forms.
Therefore it suffices to prove an appropriate set of abstract `Kummer congruences' for the Fourier coefficients 
\[         C \big( \aaa, \Phi( \g_{\Ind\chi}, r-k+1) \big| U(\mm'\mm_0^{-1}) \big)       \] 
as the character $\chi$ varies. 
\par
Applying the formula from the end of \S\ref{section3}, we have 
\[   
 C\big(\aaa,  \Phi(\g_\rho, r-k+1)\big| U(\mm'\mm_0^{-1})\big)   \;\equiv\;    N(\tilde{t}_\lambda)^{ -\frac{1}{2}-r }  (-1)^{d(r-k+1)}    \;\;\; \times
\]
\[        
 \sum_{\xi = \xi_1 + \xi_2} b_\lambda(\xi_1) \; N(\xi_1)^{r-k+1}   \sum_{\tilde{\xi_2} = \tilde{b}\tilde{c}} \mathrm{sign}( N(\tilde{b}) )^{k}    \;      
 N(\tilde{b})^{2r-k} \; (\eta \omega^{-1})(\tilde{c}) \; \mod N(\mm'\mm_0^{-1})    
\]
where $\aaa = \xi \tilde{t}_\lambda^{-1}$. All the terms in this sum are $p$-integral, and the only ones which depend on $\chi$ are the Fourier coefficient $b_\lambda(\xi_1)$ and the value of the Hecke character $\omega$. By definition of the theta-lift modular form $\g_{\Ind \chi}$, we may write
\[ b_\lambda (\xi_1)  N(\tilde{t}_\lambda)^{-1/2}  \;=\;  \sum_\mathfrak{A} \chi(\mathfrak{A}) .  \] 
where the sum ranges over all ideals $\mathfrak{A}$ of $\OO_K$ such that $N_{K/F}(\mathfrak{A}) = \xi_1  N(\tilde{t}_\lambda)^{-1}$.
Further, the Hecke character $\omega$ is given by
\[    \omega(\aaa)   \; = \; \theta_{K/F}(\aaa) \; \chi (\aaa \OO_K)    \]
where $\theta_{K/F}$ is the quadratic Hecke character defined by
\begin{equation*}
 \theta_{K/F} (\qq) = 
\begin{cases}
1&\text{if $\qq$ splits in $K/F$}\\
-1&\text{if $\qq$ is inert in $K/F$}\\
0&\text{if $\qq$ ramifies in $K/F$} 
\end{cases}
\end{equation*}
for any prime ideal $\mathfrak{q}$ (these properties are easily verified from the discussion of theta-lifts in \cite[\S5]{Shimura}). Therefore we have
\[   
 C \big( \aaa,  \Phi(\g_\rho, r-k+1) \big| U(\mm'\mm_0^{-1}) \big)   \;\equiv\;   \sum_{\xi_1, b, c, \mathfrak{A}}  w(\xi_1,b,c) \; \chi(\mathfrak{A}) \; \chi^{-1}( \tilde{c} \, \OO_K)   \!\!\!\!  \mod N(\mm'\mm_0^{-1})   
\]
for some algebraic numbers $w(\xi_1,b,c)$ which are all $p$-integral.
\par
Suppose that we have coefficients $B_\chi \in \C_p$ (almost all zero) satisfying
\[
\sum_{\chi} B_\chi \;  \chi(\sigma) \; \in p^m \,\OCp
\]
for all $\sigma \in \Gcal(\MM_0)$. For the finite set of characters $\chi$ such that $B_\chi \neq 0$, we choose $\mm'$ large enough that $p^m \cc(\g_{\Ind \chi})\mm_0^2 | \mm'$ for each $\chi$. 
We also fix ideals $\CC_c$ such that $\chi(\CC_c) = \chi^{-1}(\tilde{c} \, \OO_K)$ for each $\chi$; it is clear this can be done, since the characters are all defined modulo $\MM_0 \, p^R$ for $R$ sufficiently large. Then we may write
\[
\sum_{\chi} B_\chi  \; C \big( \aaa,  \Phi(\Ind \chi, r-k+1) \big| U(\mm'\mm_0^{-1}) \big) 
\]
\[
 \equiv \;\,  \sum_{\chi} \sum_{\xi_1, b, c, \mathfrak{A}} B_\chi \; w(\xi_1,b,c) \; \chi(\mathfrak{A}\, \CC_c) \;   \mod N(\mm'\mm_0^{-1})   .
\]
By assumption, $\sum_{\chi} B_\chi \, \chi(\mathfrak{A}\, \CC_c) \equiv 0 \mod p^m$; and we chose $\mm'$ divisible by $p^m$, so we conclude that the whole sum lies in $p^m \,\OCp$.
Thus, we have checked the Kummer congruences for these special values, and we conclude that they define a bounded $p$-adic measure on $\Gcal(\MM_0)$.
\section{Integrality of Special Values}\label{section5}
In order to establish our set of congruences, we need to know that our $p$-adic measures take integral values. In fact, it appears that they become increasingly $p$-integral with the discriminant of the totally real field $F$. We now impose our hypothesis on the congruence properties of $\f$:
\\\\
{\it Hypothesis (Cong):} there exists no congruence modulo $\MCp$ between $\f$ and another Hilbert modular form which lies outside the $\f$-isotypic component of $\mathcal{M}_k(\cc(\f),\eta)$.
\\\\
Here, when we refer to a congruence of Hilbert modular forms, we mean a congruence of their Fourier expansions. We will discuss why this hypothesis is necessary at the end of \S\ref{section6}.
\begin{prop}\label{integralityprop}
If $p>k-2$ and $\f$ satisfies Hypothesis (Cong), then we have
\[
   \nu_p \left(  \epsilon_F(\rho, 1-r) \; \frac{\Psi(\f,\g_\rho^\iota,r)}{\left< \f, \f \right>} \right) \; \geq \; \nu_p\left(D_F^{k-2} \right)
\]
for each integer $1 \leq r \leq k-1$, where $\nu_p$ is the $p$-adic valuation and $D_F$ is the discriminant of $F$.
\end{prop}
{\it Proof.} 
We have the integral representation
\[
\Psi(\f,\g^\iota,s)  \;\; = \;\;  (-1)^{d(s-k+1)}   \;   2^{dk}   \;   i^{d(1-k)}   \;  N\big( \cc(\f) \cc(\g) \dd^2 \big)^{(k-1)/2-s}  
 \;  N\big( \cc(\g) \big)^{1-k/2} 
\]
\[
\times \quad \Lambda(\f^\iota) \; \Big< \f , \Theta(s-k+1)  \Big>_{\cc(\f)}   
\]
where
\[
\Theta(s) \; = \; \mathcal{H}ol \Big( \left(\g^\iota \big| J_{\cc(\f)\cc(\g)} \right) . \, E_{k-1}(s,\eta\omega^{-1}) \Big) \Big| U\big(\cc(\g)\big) .
\]
This follows after substituting $\FF=\f$ and $\G=\g^\iota$ into Theorem \ref{rep1} and applying the formula
\begin{eqnarray*}
\big<  \f^\iota, \Theta | J_{\cc(\f)} \big>    &=&   \big<  \f^\iota| J_{\cc(\f)} , \Theta \big>
\\
    &=&  \Lambda(\f^\iota) \big<  \f , \Theta \big> .
\end{eqnarray*}
The identity above holds because $J_\cc$ is self-adjoint with respect to the Petersson inner product, and $\f^\iota|J_{\cc(\f)} = \Lambda(\f^\iota)\, \f$ for any primitive $\f$.
\par
We know $\Lambda(\f^\iota)$ is a root of unity, and we always assume that $\p \nmid 2\cc(\f)$; therefore we can write
\[
     \frac{\Psi(\f,\g_\rho^\iota,r)}{\left< \f, \f \right>}   \; = \;  \;\unit \; \times  N\big( \dd^2 \big)^{(k-1)/2-r}  
   N\big( \cc(\g) \big)^{1/2-r} 
 \frac{ \big< \f , \Theta(r-k+1) \big>_{\cc(\f)} }{\left< \f, \f \right>} 
\]
for $1 \leq r \leq k-1$. Recall that
\[
 \epsilon_F(\rho,1-r) \; = \; i^d\; \Lambda(\g_\rho) \; N\big( \cc(\g_\rho) \dd^2 \big)^{r-1/2}  .
\]
Also, the norm of the different $\dd$ is equal to the absolute discriminant $|D_F|$, so
\[
 \epsilon_F(\rho,1-r)  \; \frac{\Psi(\f,\g_\rho^\iota,r)}{\left< \f, \f \right>} \; = \; \unit \;\times\; D_F^{k-2} \; \frac{ \big< \f , \Theta(r-k+1) \big>_{\cc(\f)} }{\left< \f, \f \right>_{\cc(\f)}}  .
\]
By choosing a basis for the finite-dimensional vector space $\mathcal{M}_k(\cc(\f),\eta)$ which includes $\f$, we may write
\[
         \Theta(r-k+1)  \; = \; c \, \f \; + \; \sum_{\f_i \neq \f} c_i \; \f_i \big| \mathfrak{b}_i    
\]
for scalars $c$ and $c_i \in \Qbar$ (almost all zero), and primitive forms $\f_i \in \mathcal{M}_k(\aaa_i,\eta)$ such that 
$\aaa_i \bb_i$ divides $\cc(\f)$. We deduce that
\[    \frac{ \big< \f , \Theta(r-k+1) \big>_{\cc(\f)} }{\left< \f, \f \right>_{\cc(\f)}}  
\; = \;  c + \sum_{\f_i \neq \f} c_i \,  
 \frac{\left< \f , \f_i \big| \bb_i   \right>_{\cc(\f)}}{\left< \f, \f \right>_{\cc(\f)}}  .  \]
For each $i$ we have $\left< \f  , \f_i   \right>_{\cc(\f)} = 0$ as $\f$ and $\f_i$ are distinct primitive forms. 
This implies that $\left< \f  , \f_i \big| \bb_i \right>_{\cc(\f)}=0$ (one can see this from \cite[proposition 4.13]{Shimura}) so we have
\[
          \frac{ \big< \f , \Theta(r-k+1) \big>_{\cc(\f)} }{\left< \f, \f \right>_{\cc(\f)}} \; = \;   c . 
\]
Therefore we obtain the equation
\[
          \epsilon_F(\rho,1-r)  \; \frac{\Psi(\f,\g_\rho^\iota,r)}{\left< \f, \f \right>} \; = \; \unit \;\times\; D_F^{k-2} \; c 
\]
and to prove the proposition it suffices to show that $c$ is $p$-integral.
\par
Suppose not; then $c^{-1} \equiv 0 \mod \MCp$. We have an explicit formula for the Fourier coefficients of 
\[
\Theta(r-k+1) \; = \; \mathcal{H}ol \Big( \left(\g^\iota \big| J_{\cc(\f)\cc(\g)} \right) \cdot \, E_{k-1}(r-k+1,\eta\omega^{-1}) \Big) \, \Big|\, U\big(\cc(\g)\big) 
\]
almost identical to that for $C(\mm, \Phi(\g, r-k+1))$ in Section \ref{section3}; they can be written as $p$-integral linear combinations of the polynomials $P_{r-k+1}(\xi_{2,\nu}, \xi_{\nu})$. 
Since we assume $p>k-2$, these polynomials have $p$-integral coefficients, so it is clear that $\Theta(r-k+1)$ has a $p$-integral Fourier expansion. Therefore 
\[
      c^{-1} \; \Theta(r-k+1) \; \equiv \; 0 \mod \MCp .
\]
This means that
\begin{eqnarray*}
 \f   &=&  c^{-1} \, \Theta(r-k+1) \; - \; \sum_{\f_i\neq \f} c^{-1}  c_i \; \f_i \big| \mathfrak{b}_i   
\\ 
&\equiv&  - \sum_{ \f_i\neq \f} c^{-1}  c_i \; \f_i \big| \mathfrak{b}_i \; \quad  \mod \MCp  ,
\end{eqnarray*}
and we see that $\f$ is congruent modulo $\MCp$ to an element of $\mathcal{M}_{k}(\cc(\f),\eta)$ which does not lie in the $\f$-isotypic component. This contradicts Hypothesis (Cong) and we have proved the proposition.
\section{False Tate Curve Congruences}\label{section6}
We now restrict our attention to the setting of the introduction. We consider the $p$-power cyclotomic field $K_n = \Q(\mu_{p^n})$ and its maximal real subfield $F_n = K_n^+$, for an integer $n \geq 1$. The field $F_n$ is totally real, and from now on we assume our primitive cusp form $\f \in \mathcal{M}_k(\cc(\f),\eta)$ is defined over $F_n$. The prime $p$ is totally ramified in the extension $F_n/\Q$, and we set $\p$ to be the unique prime ideal of $\OO_{F_n}$ above $p$.
\par
We fix a false Tate curve extension of $\Q$, defined by
\[
\QFT \; := \; \bigcup_{n \geq 1} \Q \big( \mu_{p^n},  \sqrt[p^n]{\Delta} \big)
\]
where $\Delta$ is a $p$-power free integer.
Recall from \S1 that $\GFT := \Gal(\QFT/\Q)$ has a unique irreducible self-dual representation $\rho_{n,\Q} =\Ind_{K_n}^\Q \chi_n$ of dimension $p^n-p^{n-1}$, for each $n\geq 1$, and all irreducible representations of $\GFT$ may be written in the form $\rho_{n,\Q} \otimes \psi$ for some $n\geq 0$ and some finite-order character $\psi$ of $U^{(n)}$. 
\par
As before, we have a theta-lift Hilbert modular form over $F_n$, corresponding to the Hecke character $\chi_n$. We denote this by $\g_{\rho_n}$, identifying it with the two-dimensional Artin representation $\rho_n =\Ind_{K_n}^{F_n} \chi_n$. We point out that this notation is canonical, as $L(\g_{\rho_n},s)$ is equivalent to the Artin $L$-series $L(\rho_n/F_n,s)$ and we have $L(\rho_n/F_n,s)= L(\rho_{n,\Q},s)$ by the Artin formalism.
\par
Additionally, for $0 \leq j \leq n$, we have the restricted character $\chi_{j,n} := \Res_{K_n}(\chi_j)$ and a corresponding induced representation $\rho_{j,n} := \Ind^{F_n}_{K_n} \chi_{j,n}$. In fact, the theta-lift modular form $\g_{\rho_{j,n}}$ coincides with the automorphic base-change of $\g_{\rho_j}$ from $F_j$ to $F_n$. 
\par
Before we state the main result of this section, we must comment on our conventions for local epsilon factors: we follow those of Deligne from \cite{Deligne} (see also Tate \cite{Tate}). Recall that we may write the global epsilon factor as a product: 
\[
\epsilon_F(\rho,s)  \;=\; \prod_v \epsilon_F(\rho,s)_v
\]
where $v$ ranges over all places of $F$ and $\epsilon_F(\rho,s)_v$ is a local epsilon factor at $v$. In fact, the local factor at $v$ depends on a choice of additive character of $F_v$ and a Haar measure on $F_v$. We assume we have fixed these as in \cite[(3.5)]{Tate} so that the above product formula holds.
In particular, at the archimedean places we can use the standard characters and measures given by \cite[(3.2.4) and (3.2.5)]{Tate}.
\begin{prop}\label{Katolemma}
Suppose that $p>k-2$ and that $\f$ and satisfies Hypotheses (Cong) and (Ord). Suppose further that $\f$ has rational Fourier coefficients;
then for each $0 \leq j \leq n$ there exists a unique element $\Lp(\f,\rho_{j,n},r) \in \Z_p[[U^{(n)}]]$ with the property
\[
\psi \Big(  \Lp(\f,\rho_{j,n},r)  \Big) 
\;=\;  
\frac{\epsilon_{F_n}(\rho_{j,n}\otimes\psi,1-r)_\p}{ \alpha(\p)^{f(\rho_{j,n}\otimes\psi, \p)}} 
\;\times \;
\frac{ P_\p(\rho_{j,n}\otimes\psi^{-1}, \, \alpha(\p)^{-1} p^{r-1}) }{ P_\p(\rho_{j,n}\otimes\psi, \, \alpha(\p) \, p^{-r})} 
\]
\[
\;\times \;
\frac{ \Psi_S(\f, \g_{\rho_{j,n}\otimes\psi}^\iota, r) }{ D_F^{k-2} \; \left< \f,\f \right>_{\cc(\f)} } 
\]
for each character $\psi : U^{(n)}\rightarrow \C^\times$. Here we have written $\epsilon_{F}(\rho,1-r)_\p$ for the local epsilon factor at $\p$ (normalised as above). Also, we write $P_\p(\rho,T)$ for the local polynomial of $\rho$ at $\p$, $f(\rho,\p)$ for the power of $\p$ dividing the conductor of $\rho$, and  $S$ for the set of all primes dividing $p \Delta \OO_{F_n}$.
\par
Further, we have a congruence
\[
 \Lp(\f,\rho_{j,n},r) \; \equiv\;  \Lp(\f,\rho_n,r) \mod p\,\Z_p[[U^{(n)}]]
\]
for each $0 \leq j \leq n$.
\end{prop}
{\it Proof.} 
We will prove this lemma using the measure $\mu_{\mm_0}(\f,r)$ which we constructed in Section \ref{section4}. We observe that 
\[
    \rho_{j,n} \otimes \psi \; = \; \Ind_{K_n}^{F_n}(\chi_{j,n} \cdot \psi).
\] 
The characters $\chi_{j,n}$ and $\psi$ extend naturally to $\Gcal(\MM_0)$, where $\MM_0$ is the product of $\PP$ and all primes dividing $\Delta\OO_{K_n}$. 
Therefore we may define a bounded $p$-adic measure $\mu(\f,\rho_{j,n},r)$ on $U^{(n)}$ by setting
\[
\int_{ U^{(n)}} \psi \; d\mu(\f,\rho_{j,n},r)  \;=\;  \frac{1}{ D_F^{k-2}} \, \int_{\Gcal(\MM_0)} \chi_{j,n} \cdot \psi \; d\mu_{\mm_0}(\f,r) 
\]
for each finite-order character $\psi : U^{(n)}\rightarrow \C^\times$.
Under our hypotheses we may apply Proposition \ref{integralityprop} to see that these values are integral. This implies the existence of an element $\widetilde{\Lp}(\f,\rho_{j,n},r) \in \OCp[[U^{(n)}]]$ which produces the value above when evaluated at the character $\psi$ (making the usual identification of integral $p$-adic measures with elements of the Iwasawa algebra).
\par
Comparing the values of the measure $\mu_{\mm_0}(\f,r)$ with those proposed for $\Lp(\f,\rho_{j,n},r)$, we see that our desired element is equal to
\[
\Lp(\f,\rho_{j,n},r) \;=\; \gamma_{j,n}(\f,r) \,\cdot\, \widetilde{\Lp}(\f,\rho_{j,n},r)
\]
where $\gamma_{j,n}(\f,r) \in \OCp[[U^{(n)}]]$ is defined by the interpolation property
\[
 \psi \big( \gamma_{j,n}(\f,r) \big)  \;=\;
 \prod_{\qq | \Delta\OO_{F_n}} \frac{\alpha(\qq)^{\ord_\qq (\cc(\rho_{j,n}))}}{ \epsilon_{F_n}(\rho_{j,n}\otimes\psi,1-r)_\qq }.
\]
Since $p$ and $\Delta$ are coprime, each $\alpha(\qq)$ is a $p$-adic unit. Further, since $\psi$ is unramified at $\qq | \Delta.\OO_{F_n}$, we may apply \cite[3.4.6]{Tate} to write
\[
\epsilon_{F_n}(\rho_{j,n}\otimes\psi,1-r)_v \; = \; \psi(\qq^{A_{j,n}}) \,\cdot\, \epsilon_{F_n}(\rho_{j,n},1-r)_v
\]
where the exponent $A_{j,n}$ does not depend on $\psi$. It is clear that such a $\gamma_{j,n}(\f,r)$ exists, which establishes the existence of $\Lp(\f,\rho_{j,n},r) \in \OCp[[U^{(n)}]]$.
\par
To demonstrate that this element actually lies in $\Z_p[[U^{(n)}]]$, we apply \cite[theorem 4.2]{Shimura} which shows that the value
\[
 \tau(\omega)^{-1} \;
\frac{ \Psi_S(\f, \g, r) }{  \left< \f,\f \right>_{\cc(\f)} } 
\]
is Galois-equivariant with respect to $(\f,\g)\mapsto (\f^\sigma,\g^\sigma)$ for all $\sigma \in \mathrm{Aut}(\C)$; here, $\tau(\omega)$ denotes the Gauss sum of the Hecke character $\omega$ associated to $\g$. In the case $\g=\g_{\rho_{j,n}\otimes\psi}^\iota$, it is clear that $\tau(\omega)^{-1}$ will be a rational multiple of the epsilon factor $\epsilon_{F_n}(\rho_{j,n}\otimes\psi,1-r)$ and we have the same Galois-equivariance property for
\[
\epsilon_{F_n}(\rho_{j,n}\otimes\psi,1-r) \;
\frac{ \Psi_S(\f, \g_{\rho_{j,n}\otimes\psi}^\iota, r) }{  \left< \f,\f \right>_{\cc(\f)} } .
\]
Further, we can replace the global epsilon factor by its local counterpart at $\p$ without affecting this property. We may do this because the conductors of the representations $\rho_{j,n}$ away from $p$ are squares, and the representations are self-dual so their local root numbers are $\pm 1$. Therefore, the local factors at primes other than $\p$ lie in $\Q$ and they do not affect the rationality of the values.
\par
Since we now assume that $\f$ has rational Fourier coefficients, and $\rho_{j,n}$ may be realised over $\Q$, this Galois-equivariance property shows that the values of $\Lp(\f,\rho_{j,n},r)$ lie in $\Z_p$.
\par
It remains to prove that these elements are congruent modulo $p$ as $j$ varies. The character $\chi_{j,n}$ takes values in $\mu_{p^n}$ for all $0 \leq j \leq n$. Therefore we have $\chi_{j,n} \equiv \chi_n \mod \MCp$, which implies
\[
\frac{1}{ D_F^{k-2}} \int_{\Gcal(\MM_0)} \psi \,.\, \chi_{j,n} \; d\mu_{\mm_0}(\f,r) \; \equiv \; \frac{1}{ D_F^{k-2}} \int_{\Gcal(\MM_0)} \psi \,.\, \chi_{n}\; d\mu(\f,\rho_n,r) \mod \MCp
\]
for any $\psi$ (as this measure is integral). So the values $\psi \big( \Lptilde(\f,\rho_{j,n},r) \big)$ and 
$\psi \big( \Lptilde(\f,\rho_{n},r) \big)$ are congruent for all $\psi$, implying the elements of $\OCp[[U^{(n)}]]$ themselves are congruent. Similarly, we have a congruence
\[
\gamma_{0,n}(\f,r) \; \equiv \; \gamma_{1,n}(\f,r) \; \equiv \,...\, \equiv \; \gamma_{n,n}(\f,r) \mod \MCp[[U^{(n)}]];
\]
this fact is analogous to Claim ($\star$) in \cite{DelbourgoWard}, and follows easily from the proof of that assertion. This implies that
$\Lp(\f,\rho_{j,n},r)$ and $\Lp(\f,\rho_n,r)$ are congruent modulo $\MCp \cap \Z_p = p\, \Z_p$.
\par
As an immediate consequence of this proposition, we have proved Theorem \ref{theorem1}.
Now, let $\varphi$ and $N_{i,j}$ be the maps defined in \S\ref{Intro}, and put $a_j= \Lp(\f,\rho_j,r)$ for a critical value $1 \leq r \leq k-1$. A simple induction argument (given in detail in \cite{DelbourgoWard}) allows us to deduce the following result from Proposition \ref{Katolemma} (which establishes Theorem \ref{theorem2}).
\begin{corollary}\label{Katocorol}
Assume that $p>k-2$, $\f$ has rational Fourier coefficients, and $f$ satisfies Hypotheses (Cong) and (Ord). Then we have the `false Tate curve' congruence
\[     
\prod_{1\leq j \leq n} N_{j,n}\left(  \frac{a_j}{N_{0,j}(a_0)}. \frac{\varphi \circ N_{0,j-1}(a_0)}{\varphi(a_{j-1})} \right)^{p^j}  \, \equiv \, 1 \mod p^{n+1}  
\]
for each $n\geq 0$.
\end{corollary}
\noindent
{\bf Remark:}
It is not difficult to find cases in which Hypothesis (Cong) fails and Theorem \ref{theorem2} does not hold.  
Let us give an explicit example: taking $p=3$ we have $F_1=\Q$, so if we work over this field we are reduced to the case of elliptic modular forms. This allows us to compute the values of $\Lp(f,\rho_{j,1},r)$ using MAGMA \cite{MAGMA} (we will discuss the methods used in \S\ref{section7}). 
\par
We have a primitive cusp form $f$ in $S_4^{\mathrm{new}}(\Gamma_0(19))$, with $q$-expansion $f(z)=q - 3q^2 - 5q^3 + q^4 +...$. This cusp form is congruent modulo $3$ to another modular form $\tilde{f}(z) = q + 9q^2 + 28q^3 + 73q^4 +...$ at level $19$. 
We compute the value of $\Lp(f,\rho_{0,1},r)\big|_{r=1}$ and  $\Lp(f,\rho_{1,1},r)\big|_{r=1}$ evaluated at the trivial character:
\[
\mathbf{1}\Big(\Lp(f,\rho_{0,1},1)\Big) \,=\, 3^2+2.3^5 + O(3^{7})
, \;\;\;
\mathbf{1}\Big(\Lp(f,\rho_{1,1},1)\Big) \,=\, 1+2.3^1 + 2.3^2 + O(3^4).
\]
If Theorem \ref{theorem2} held here, these values would be congruent modulo $3$, which is not the case; the congruence at $r=3$ also fails.
\par
Let us say briefly why this is not surprising. The complex period in the interpolation formula for $\Lp(\f,\rho_{j,1},r)$ is the Petersson inner product $\left< \f,\f \right>_{\cc(\f)}$, which is equal to a twisted adjoint $L$-value of $\f$, up to certain simple factors (see \cite[theorem 7.1]{Hida-Tilouine}). 
\par
It is known that the adjoint $L$-series is closely related to the congruence module of the Hecke eigenform $\f$; this is discussed at length by Doi, Hida and Ishii in \cite{Doi-Hida-Ishii}. In particular, $p$ will be a congruence prime for $\f$ if and only if $p$ divides the algebraic part of the adjoint $L$-value. This relationship is conjectural for Hilbert modular forms in full generality, but has been proved under certain assumptions: see \cite{Ghate} and \cite{Dimitrov} for example.
\par
As a consequence, if $\p$ is a congruence prime for $\f$ then $\p$ should divide the algebraic part of the adjoint $L$-value at $s=1$, making the values of $\Lp(\f,\rho_{j,n},r)$ less $p$-integral. Therefore Hypothesis (Cong) may be a natural condition to impose, when we use this particular automorphic period.
\section{Numerical Examples over $\Q$}\label{section7}
We have established Proposition \ref{integralityprop} and Corollary \ref{Katocorol} subject to the assumption $p>k-2$.
More precisely: to prove Proposition \ref{integralityprop} we needed the coefficients of the polynomials $P_{r-k+1}(\xi_{2,\nu}, \xi_{\nu})$ to be $p$-integral, and their denominators could only be divisors of $(k-2)!$. 
Despite the fact that we were unable to remove this hypothesis, we expect that our results should hold even when $p\leq k-2$ and we will conclude our paper by giving some numerical evidence for this.
\par
Unfortunately it is difficult to compute convolution $L$-values of Hilbert modular forms in general, but if we restrict to the case $p=3$ and $n=1$, our base field is $F_1=\Q$ and we can work with elliptic modular forms. 
Using the $L$-series functions in the computer package MAGMA \cite{MAGMA}, we calculate the value of $\psi\big(\Lp(f,\rho_{j,1},r)\big)$ in the case $\psi=\mathbf{1}$ for several examples of $f$ (where we write $f$ instead of $\f$ to emphasise the fact that we are now working with a classical modular form). To compute the interpolation factors, we use the methods of T. and V. Dokchitser which are described in the paper \cite{VladTim}.
\par
To compute the Petersson inner product $\left< f,f \right>_{N}$ for a cusp form $f \in S_k(\Gamma_0(N))$ we use a well-known formula which relates it to the adjoint $L$-series:
\[
\left< f,f \right>_{N} \; = \; 2^{-2k} \; (k-1)! \; N \; \phi(N) \; \frac{L \big( \Ad (f),1 \big)}{\pi^{k+1}} .
\] 
We can then compute the algebraic parts of the $L$-values for $\rho= \rho_{0,1}$ and $\rho_{1,1}$, which we write as follows:
\[
\Psi^*_f(\rho,r)   \;=\; \sqrt{N_\rho} \cdot \frac{ \Psi_S(f, g_{\rho}, r) }{ \left< f,f \right>_{N} } .
\]
Here $N_\rho$ denotes the conductor of $\rho$. Finally we compute the values
\[
\Lcal_f(\rho,r) = 
\frac{\epsilon(\rho,1-r)_p}{ \alpha(p)^{f(\rho, p)}} 
\;\times \;
\frac{ P_p(\rho^\vee, \, \alpha(p)^{-1} p^{r-1}) }{ P_p(\rho, \, \alpha(p) \, p^{-r})} 
\;\times \;
\frac{ \Psi_S(f, g_{\rho}, r) }{ \left< f,f \right>_{N} } .
\]
This is equal to the evaluation of $\Lp(f,\rho,r)$ at the trivial character. By Proposition \ref{Katolemma} we expect the congruence
\[
\Lcal_f(\rho_{0,1},r) \;\equiv\; \Lcal_f(\rho_{1,1},r) \mod 3 .
\]
We tested two primitive cusp forms of weight $6$ for which Hypothesis (Cong) is satisfied. 
Table \ref{table6N5} shows our results for the newform $f \in S_6^{\mathrm{new}}(\Gamma_0(5))$ which has $q$-expansion $q + 2q^2 - 4q^3 - 28q^4 +...$ and Table \ref{table6N7} shows the same data for $f \in S_6^{\mathrm{new}}(\Gamma_0(7))$ having $q$-expansion $q - 10q^2 - 14q^3 + 68q^4+...$. 
As $p=3$ the assumption $p>k-2$ is not satisfied, but we still observe the congruence for each critical value of $r$. 
\begin{table}[ht]
\caption{values for $f$ of weight $6$, level $5$, with $\Delta = 2$}
\centering
\begin{tabular}{cccccc}
\hline  
$r$ & $\Psi^*_f(\rho_{0,1},r)$ & $\Psi^*_f(\rho_{1,1},r)$  &  $\Lcal_f(\rho_{0,1},r)$ & $\Lcal_f(\rho_{1,1},r)$     
\\[1ex]
\hline 
 1 &  $\frac{2^{5} . 5^{5} }{31^{1}}$  & $2^{14} . 5^{5} . 661^{1}$  & $1.3^{0} + 1.3^{1} +  O(3^{3})$   & $1.3^{0} + 2.3^{1} +  O(3^{2})$
 \\[1ex]
 2 & $\frac{2^{2} . 5^{2} }{ 3^{1}}$  & $\frac{2^{5} . 5^{3} . 1759^{1}}{  3^{3}}$  & $2.3^{0} + 2.3^{2} +  O(3^{3})$   & $2.3^{0} + 2.3^{1} + O(3^{2})$
 \\[1ex]
 3 & $\frac{2^{3} }{ 3^{2}}$  & $\frac{2^{5} . 5^{2} }{ 3^{6}}$  & $1.3^{0} + 2.3^{2} + O(3^{4})$   & $1.3^{0} + 1.3^{1} + O(3^{2})$
 \\[1ex]
 4 & $\frac{2^{2} }{ 3^{3}}$  & $\frac{2^{1} . 5^{1} . 1759^{1} }{ 3^{9}}$  & $2.3^{0} + 1.3^{1} + O(3^{2})$   & $2.3^{0} + 2.3^{4} + O(3^{5})$
 \\[1ex]
 5 & $\frac{2^{5} . 5^{1} }{ 3^{4} . 31^{1}}$  & $\frac{2^{6} . 5^{1} . 661^{1} }{ 3^{12}}$  & $1.3^{0} + 2.3^{1} + O(3^{2})$   & $1.3^{0} + 2.3^{2} + O(3^{3})$
 \\[1ex]
\hline
\end{tabular}
\label{table6N5}
\end{table}
\begin{table}[ht]
\caption{values for $f$ of weight $6$, level $7$, with $\Delta = 2$}
\centering
\begin{tabular}{cccccc}
\hline 
$r$ & $\Psi^*_f(\rho_{0,1},r)$ & $\Psi^*_f(\rho_{1,1},r)$  &  $\Lcal_f(\rho_{0,1},r)$ & $\Lcal_f(\rho_{1,1},r)$     
\\[1ex]
\hline 
 1 & $\frac{2^{8} . 7^{2} . 19^{1} }{ 3^{1} . 43^{1}}$  & $2^{10} . 5^{2} . 7^{3} . 13^{3}$  & $2.3^{0} + 1.3^{1} + O(3^{2})$   & $2.3^{0} + 1.3^{2} + O(3^{3})$
 \\[1ex]
 2 & $\frac{2^{2} . 7^{1} }{ 3^{1}}$  & $\frac{2^{8} . 7^{2} . 181^{1} }{ 3^{3}}$  & $2.3^{0} + 1.3^{3} + O(3^{4})$   & $2.3^{0} + 2.3^{1} + O(3^{2})$
 \\[1ex]
 3 &  $0$ &  $0$ & $0$   & $0$
 \\[1ex]
 4 & $\frac{2^{2} }{ 3^{3} . 7^{1}}$  & $\frac{2^{4} . 181^{1} }{ 3^{9}}$  & $1.3^{0} + 2.3^{1} + O(3^{3})$   & $1.3^{0} + 1.3^{4} + O(3^{5})$
 \\[1ex]
 5 & $\frac{2^{8} . 19^{1} }{ 3^{5} . 7^{2} . 43^{1}}$ & $\frac{2^{2} . 5^{2} . 13^{3} }{ 3^{12} . 7^{1}}$ & $1.3^{0} +   1.3^{1} +  O(3^{2})$ & $1.3^{0} +   2.3^{1} +  O(3^{2})$
 \\[1ex]
\hline
\end{tabular}
\label{table6N7}
\end{table}
\clearpage
\bibliographystyle{siam}
\bibliography{REFERENCES}
\end{document}